\documentclass[12pt]{article}
\usepackage{amssymb,amsmath}
\def\d{\delta}
\def\D{\Delta}
\def\UU{{\mathcal U}}
\def\WW{{\mathcal W}}
\def\LL{{\mathcal L}}

\def\Der{{\rm Der}}
\def\Inn{{\rm Inn}}

\def\Ker{{\rm Ker}}

\def\Im{{\rm Im}}

\def\cl{\centerline}

\def\rar{\rightarrow}

\def\vs{\vspace*}

\def\ni{\noindent}
\def\VV{\mathcal {V}}

\def\WW{\mathcal {W}}

\def\Z{\mathbb{Z}}

\def\C{\mathbb{C}}
\def\QED{\hfill$\Box$}

\def\QED{\hfill$\Box$}
\numberwithin{equation}{section}
\newtheorem{theo}{Theorem}[section]
\newtheorem{defi}[theo]{Definition}

\newtheorem{coro}[theo]{Corollary}
\newtheorem{lemm}[theo]{Lemma}
\newtheorem{prop}[theo]{Proposition}
\newtheorem{clai}{Claim}

\def\adddot{$\!\!\!${\bf.}\ \ }

\begin{document}
\cl{{\Large \bf Lie bialgebra structures on the $W$-algebra
$W(2,2)$} \footnote {Supported by NSF grants 10471091, 10671027 of
China, ``One Hundred Talents Program'' from University of Science
and Technology of China.\\[2pt] \indent Corresponding E-mail:
sd\_junbo@163.com} } \vs{6pt}

\cl{Junbo Li$^{*,\dag)}$, Yucai Su$^{\ddag)}$} \cl{\small
$^{\dag)}$Department of Mathematics, Changshu Institute of
Technology, Changshu 215500, China}  \cl{\small
$^{\ddag)}$Department of Mathematics, University of Science and
Technology of China Hefei 230026, China} \cl{\small E-mail:
sd\_junbo@163.com, ycsu@ustc.edu.cn} \vs{6pt}

{\small
\parskip .005 truein
\baselineskip 3pt \lineskip 3pt

\noindent{{\bf Abstract.} Verma modules over the $W$-algebra
$W(2,2)$ were considered by Zhang and Dong, while the Harish-Chandra
modules and irreducible weight modules over the same algebra were
classified by Liu and Zhu etc. In the present paper we shall
investigate the Lie bialgebra structures on the referred algebra,
which are shown to be triangular coboundary.\vs{5pt}

\noindent{\bf Key words:} Lie bialgebras, Yang-Baxter equation, $W$
algebra $W(2,2)$.}

\noindent{\it Mathematics Subject Classification (2000):} 17B05,
17B37, 17B62, 17B68.}
\parskip .001 truein\baselineskip 6pt \lineskip 6pt

\vs{18pt}

\cl{\bf\S1. \
Introduction}\setcounter{section}{1}\setcounter{equation}{0} It is
well known that the notion of Lie bialgebras was originally
introduced by Drinfeld in 1983 (cf.~\cite{D1}) during the search for
the solutions of the Yang-Baxter quantum equation. During the resent
rears, there have appeared several papers on Lie bialgebras (e.g.,
\cite{LSX} and \cite{M1}--\cite{WSS}). Witt and Virasoro type Lie
bialgebras were introduced in \cite{T}, of which type Lie bialgebras
were further classified in \cite{NT}, while the generalized case was
considered in \cite{SS}. Lie bialgebra structure on generalized
Virasoro-like and Block  Lie algebras were investigated in
\cite{WSS} and \cite{LSX} respectively.

In this paper we shall investigate Lie bialgebra structures on the
$W$ algebra $W(2,2)$ introduced in \cite{ZD}, denoted by $\WW$ here,
which have been quantized in \cite{LSq} by the authors using the
method introduced in \cite{G} and generalized in \cite{HW}. This
algebra is an infinite-dimensional Lie algebra with a $\C$-basis
$\{\,L_n,\,W_n,\,c\,|\,n\in \Z\,\}$ and the following Lie brackets
(\,other components vanishing):
\begin{eqnarray}
&&[L_m,L_n]=(m-n)L_{m+n}+\frac{m^3-m}{12}\delta_{m+n,0}c,\label{BLie01}\\
&&[L_m,W_n]=(m-n)W_{m+n}+\frac{m^3-m}{12}\delta_{m+n,0}c.\label{BLie02}
\end{eqnarray}
The verma modules on $\WW$ were investigated in \cite{ZD}. Later all
irreducible weight modules with finite dimensional weight spaces and
all indecomposable modules with less than one dimensional weight
space on $\WW$ were classified in \cite{LZ}. Meanwhile, irreducible
weight modules possessing at least one nontrivial finite-dimensional
weight space were also classified in \cite{LGZ}.

Let us recall the definitions related to Lie bialgebras. For
convenience, we introduce them being assort to the notation $\LL$.
Let $\LL$ be any vector space over the complex field $\C$ of
characteristic zero. Denote by $\xi$ the {\it cyclic map} of
$\LL\otimes\LL\otimes \LL$ cyclically permuting the coordinates,
namely, $ \xi (x_{1} \otimes x_{2} \otimes x_{3}) =x_{2} \otimes
x_{3} \otimes x_{1}$ for $x_1,x_2,x_3\in\LL,$ and by $\tau$ the {\it
twist map} of $\LL\otimes\LL$, i.e., $\tau(x\otimes y)= y \otimes x$
for $x,y\in \LL.$

First one need to reformulate the definitions of a Lie algebra
and Lie coalgebra as follows.\\
\ni A {\it Lie algebra} is a pair $(\LL,\delta)$ of a vector space
$\LL$ and a linear map $\delta :\LL\otimes\LL\rar\LL$ (the {\it
bracket} of $\LL$) satisfying the conditions:
\begin{eqnarray}
\label{Lie-s-s} \!\!\!\!\!\!\!\!\!\!\!\!&&
\Ker(1-\tau) \subset \Ker\,\delta,\\
\label{Lie-j-i} \!\!\!\!\!\!\!\!\!\!\!\!&& \delta \cdot (1 \otimes
\delta ) \cdot (1 + \xi +\xi^{2}) =0 : \LL \otimes\LL\otimes\LL\rar
\LL,
\end{eqnarray}
A {\it Lie coalgebra} is a pair $(\LL,\D)$ of a vector space $\LL$
and a linear map $\D:\LL\to\LL\otimes\LL$ (the {\it cobracket} of
$\LL$) satisfying the conditions:
\begin{eqnarray}
\label{cLie-s-s} \!\!\!\!\!\!\!\!\!\!\!\!&&
\Im\,\D \subset \Im(1- \tau),\\
\label{cLie-j-i} \!\!\!\!\!\!\!\!\!\!\!\!&& (1 + \xi +\xi^{2}) \cdot
(1 \otimes \D) \cdot \D =0:\LL\to\LL\otimes\LL\otimes\LL,
\end{eqnarray}
For a Lie algebra $\LL$, we always use $[x,y]=\d(x,y)$ to denote its
Lie bracket and use the symbol ``$\cdot$'' to stand for the {\it
diagonal adjoint action}
\begin{eqnarray}\label{diag}
\mbox{$x\cdot (\sum\limits_{i} {a_{i} \otimes b_{i}}) =
\sum\limits_{i} ( {[x, a_{i}] \otimes b_{i} + a_{i} \otimes [x,
b_{i}]})\mbox{ \ \ for \ }x,a_i,b_i\in\LL.$}
\end{eqnarray}

\begin{defi}\adddot\rm
 A {\it Lie bialgebra} is a triple $(\LL,\delta,\D)$
satisfying the conditions:
\begin{eqnarray}
\label{bLie-l} \!\!\!\!\!\!\!\!\!\!\!\!&& \mbox{$(\LL, \delta)$ is a
Lie algebra},\ \
\mbox{$(\LL,\D)$ is a Lie coalgebra,}\\
\label{bLie-d} \!\!\!\!\!\!\!\!\!\!\!\!&& \mbox{$\D  \delta (x,y) =
x \cdot \D y - y \cdot \D x$ for $x, y \in\LL$ \ (compatibility
condition).}
\end{eqnarray}
\end{defi}

Denote by $\UU$ the universal enveloping algebra of $\LL$ and by $1$
the identity element of $\UU$. For any $r =\sum_{i} {a_{i} \otimes
b_{i}}\in\LL\otimes\LL$, define $r^{ij},\,c(r),\,i,j=1,2,3$ to be
elements of $\UU \otimes \UU \otimes \UU$ by (where the bracket in
(\ref{add1-}) is the commutator):
\begin{eqnarray}\!\!\!\!\!\!\!\!\!\!\!\!\!\!&\!\!\!\!\!\!\!\!\!\!\!\!\!\!&
 r^{12} = \mbox{$\sum \limits_{i}$}{a_{i} \otimes b_{i}
\otimes 1} , \ \ r^{13}= \mbox{$\sum \limits_{i}$} {a_{i} \otimes 1
\otimes b_{i}} , \ \ r^{23} = \mbox{$\sum \limits_{i}$}{1 \otimes
a_{i} \otimes
b_{i}}, \nonumber\\
\!\!\!\!\!\!\!\!\!\!\!\!\!\!&\!\!\!\!\!\!\!\!\!\!\!\!\!\!&
\label{add1-} \mbox{$c(r) = [r^{12} , r^{13}] +[r^{12} , r^{23}]
+[r^{13} , r^{23}].$}
\end{eqnarray}

\begin{defi}\adddot
\label{def2} \rm (1) A {\it coboundary Lie bialgebra} is a $4$-tuple
$(\LL, \delta, \D,r),$ where $(\LL,\delta,\D)$ is a Lie bialgebra
and $r \in \Im(1-\tau) \subset\LL\otimes\LL$ such that $\D=\D_r$ is
a {\it coboundary of $r$}, where $\D_r$ is defined by
\begin{equation}
\label{D-r}\D_r(x)=x\cdot r\mbox{\ \ for\ \ }x\in\LL.
\end{equation}

 (2) A coboundary Lie bialgebra $(\LL,\delta,\D,r)$
is called {\it triangular} if it satisfies the following {\it
classical Yang-Baxter Equation} (CYBE):
\begin{equation}
\label{CYBE} c(r)=0.
\end{equation}
\end{defi}

The main result of this paper can be formulated as follows.
\begin{theo}\adddot
\label{main} Every Lie bialgebra structure on $\WW$ is triangular
coboundary. \end{theo}

Throughout the paper, we denote by $\Z^*$ the set of all nonzero
integers, $\Z_+$ the set of all nonnegative integers and $\C^*$ the
set of all nonzero complex numbers.

\vs{12pt}

\cl{\bf\S2. \ Proof of the main results}\setcounter{section}{2}
\setcounter{theo}{0}\setcounter{equation}{0}

\vs{10pt}

The following result can be found in \cite{D1,D2,NT}.
\begin{lemm}\adddot
\label{some} Let $\LL$ be a Lie algebra and $r\in\Im(1-\tau)\subset
\LL\otimes\LL$.

 {\rm(1)} The tripple $(\LL,[\cdot,\cdot], \D_r)$ is a Lie
bialgebra if and only if $r$ satisfies CYBE $(\ref{CYBE})$.

{\rm(2)} We have
\begin{equation}
\label{add-c}(1+\xi+\xi^{2})\cdot(1\otimes\D)\cdot\D(x)=x\cdot c(r)
\mbox{\ for all\ }x\in\LL.
\end{equation}
\end{lemm}

\begin{lemm}
\label{Legr}
 \adddot Let $\WW^{\otimes n}=\WW\otimes\cdots\otimes\WW$ be
the tensor product of $n$ copies of $\WW$, and regard $\WW^{\otimes
n}$ as a $\WW$-module under the adjoint diagonal action of $\WW$.
Suppose $r\in\WW^{\otimes n}$ satisfying $x\cdot r=0$ for all
$x\in\WW$. Then $r\in\C c^{\otimes n}$.
\end{lemm} \ni{\it
Proof.~}~It can be proved directly by using the similar arguments as
those presented in the proof of Lemma 2.2 of
\cite{WSS}.\QED\vskip5pt

An element $r\in\Im(1-\tau)\subset\WW\otimes\WW$ is said to satisfy
the \textit{modified Yang-Baxter equation} (MYBE) if
\begin{equation}
\label{MYBE}x\cdot c(r)=0\ \mbox{\ for all\ }x\in\WW.
\end{equation}
As a conclusion of Lemma \ref{Legr}, one immediately obtains

\begin{coro}\label{coro1}\adddot
An element $r\in\Im(1-\tau)\subset\WW\otimes\WW$ satisfies CYBE
$(\ref{CYBE})$ if and only if it satisfies MYBE $(\ref{MYBE})$.
\end{coro}

Regard $\VV=\WW\otimes\WW$ as a $\WW$-module under the adjoint
diagonal action. Denote by $\Der(\WW,\VV)$ the set of
\textit{derivations} $D:\WW\to\VV$, namely, $D$ is a linear map
satisfying
\begin{equation}
\label{deriv} D([x,y])=x\cdot D(y)-y\cdot D(x)\mbox{ \ for \
}x,y\in\WW,
\end{equation}
and $\Inn(\WW,\VV)$ the set consisting of the derivations $v_{\rm
inn}, \, v\in\VV$, where $v_{\rm inn}$ is the \textit{inner
derivation} defined by
\begin{equation}
\label{inn} v_{\rm inn}:x\mapsto x\cdot v\mbox{ \ for \ }x\in\WW.
\end{equation}
Then it is well known that
\begin{equation} H^1(\WW,\VV)\cong\Der(\WW,\VV)/\Inn(\WW,\VV),
\end{equation}
where $H^1(\WW,\VV)$ the {\it first cohomology group} of the Lie
algebra $\WW$ with coefficients in the $\WW$-module $\VV$.

\begin{prop}\adddot
\label{proposition} $\Der(\WW,\VV)=\Inn(\WW,\VV)$, equivalently,
$H^1(\WW,\VV)=0$.
\end{prop}

\ni{\it Proof.} Note that $\WW=\oplus_{n\in\Z}\WW_n$ and
$\VV=\WW\otimes\WW=\oplus_{n\in\Z}\VV_n$ are $\Z$-graded with
\begin{equation}
\WW_n={\rm Span}\{L_n,\,W_n\,|\,n\in \Z\}\oplus\d_{n,0}\C c \mbox{ \
\ and \ } \VV_n=\mbox{$\sum\limits_{p,q\in\Z,\atop
p+q=n}$}\WW_p\otimes\WW_q\mbox{\ \,\ for \ }n\in\Z.
\end{equation}
A derivation $D\in\Der(\WW,\VV)$ is {\it homogeneous of degree
$\alpha\in\Z$} if $D(\WW_n)\subset \WW_{\alpha+n}$ for all $n\in\Z$.
 Denote $$\Der(\WW,\VV)_\alpha =\{D\in\Der(\WW,\VV)\,|\,{\rm deg\,}D=
\alpha\}\mbox{ \ \ for \ } \alpha\in\Z.$$ Let $D\in\Der(\WW,\VV)$.
For any $\alpha\in\Z$, we define the linear map
$D_\alpha:\WW\rightarrow\VV$ as follows: For any $\mu\in\WW_n$ with
$n\in\Z$, write $D(\mu)=\sum_{p\in\Z}\mu_p$ with $\mu_p\in\VV_p$,
then we set $D_\alpha(\mu)=\mu_{n+\alpha}$. Obviously, $D_\alpha\in
\Der(\WW,\VV)_\alpha$ and we have
\begin{eqnarray}\label{summable}
D=\mbox{$\sum\limits_{\alpha\in\Z}D_\alpha$},
\end{eqnarray}
which holds in the sense that for every $u\in\WW$, only finitely
many $D_\alpha(u)\neq 0,$ and $D(u)=\sum_{\alpha\in\Z}D_\alpha(u)$
(we call such a
sum in (\ref{summable}) {\it summable}).\\[8pt]

We shall prove this proposition by several claims.
\begin{clai}\adddot
\label{clai1} \rm If $\alpha\in\Z^*$, then
$D_\alpha\in\Inn(\WW,\VV)$.
\end{clai}

For $\alpha\neq 0$, denote
$\gamma=\alpha^{-1}D_{\alpha}(L_0)\in\VV_{\alpha}$. Then for any
$x_n\in\WW_{n}$, applying $D_{\alpha}$ to $[L_0,x_n]=-nx_n$, using
$D_{\alpha}(x_n)\in \VV_{n+\alpha}$, we obtain
\begin{eqnarray}\label{equa-add-1}
&&\!\!\!\!\!\!-(\alpha+n)D_{\alpha}(x_n)-x_n\cdot
D_{\alpha}(L_0)=L_0\cdot D_{\alpha}(x_n)-x_n\cdot
D_{\alpha}(L_0)=-nD_{\alpha}(x_n),
\end{eqnarray}
 i.e.,
$D_{\alpha}(x_n)=\gamma_{\rm inn}(x_n)$. Thus
$D_{\alpha}=\gamma_{\rm inn}$ is inner.

For convenience, we always use ``$\equiv$'' to denote equal modulo
$\C(c\otimes c)$ in the following.

\begin{clai}\adddot
\label{clai2} $D_0(L_0)\equiv D_0(c)\equiv 0$.
\end{clai}
\par
For any $n\in\Z$ and $x_n\in\WW_n$, applying $D_0$ to
$[L_0,x_n]=-nx_n$ and $[x_n,c]=0$ respectively, one has $x_n\cdot
D_0(L_0)=x_n\cdot c=0$. Thus by Lemma \ref{Legr}, $D_0(L_0)\equiv
D_0(c)\equiv 0$. \vskip4pt

\begin{clai}\adddot
\label{clai3}  Replacing $D_0$ by $D_0-u_{\rm inn}$ for some $u\in
\VV_0$, one can suppose $D_0(\WW)\equiv 0$.
\end{clai}
\vskip4pt
\par
For any $m\in\Z^*,\,n\in\Z$, one can write $D_0(L_m)$ and $D_0(W_n)$
as follows
\begin{eqnarray}
D_0(L_m)\!\!\!&\equiv&\!\!\!\mbox{$\sum\limits_{p\in\Z}a_{m,p}$}L_p\otimes
L_{m-p}+\mbox{$\sum\limits_{p\in\Z}b_{m,p}$}L_p\otimes W_{m-p}+a_m
L_m\otimes c+b_m c\otimes L_m\nonumber\\
&&\!\!\!+\mbox{$\sum\limits_{p\in\Z}c_{m,p}$}W_p\otimes
L_{m-p}+\mbox{$\sum\limits_{p\in\Z}d_{m,p}$}W_p\otimes W_{m-p}+c_m
W_m\otimes c+d_m c\otimes W_m,\label{0801251}\\
D_0(W_n)\!\!\!&\equiv&\!\!\!\mbox{$\sum\limits_{p\in\Z}e_{n,p}$}L_p\otimes
L_{n-p}+\mbox{$\sum\limits_{p\in\Z}f_{n,p}$}L_p\otimes W_{n-p}+e_n
L_n\otimes c+f_n c\otimes L_n\nonumber\\
&&\!\!\!+\mbox{$\sum\limits_{p\in\Z}g_{n,p}$}W_p\otimes
L_{n-p}+\mbox{$\sum\limits_{p\in\Z}h_{n,p}$}W_p\otimes W_{n-p}+g_n
W_n\otimes c+h_n c\otimes W_n,\label{0801252}
\end{eqnarray}
where all the coefficients of the tensor products are complex
numbers, and the sums are all finite. For any $p\in\Z$, the
following identities hold:
\begin{eqnarray*}
&&L_1\cdot(L_0\otimes c)=L_{1}\otimes c,\ \ \
L_1\cdot(c\otimes L_{0})=c\otimes L_{1},\\
&&L_1\cdot(W_0\otimes c)=W_{1}\otimes c,\ \ \
L_1\cdot(c\otimes W_{0})=c\otimes W_{1},\\
&&L_1\cdot(L_p\otimes L_{-p})=(1-p)L_{p+1}\otimes
L_{-p}+(1+p)L_p\otimes L_{1-p},\\
&&L_1\cdot(L_p\otimes W_{-p})=(1-p)L_{p+1}\otimes
W_{-p}+(1+p)L_p\otimes W_{1-p},\\
&&L_1\cdot(W_p\otimes L_{-p})=(1-p)W_{p+1}\otimes
L_{-p}+(1+p)W_p\otimes L_{1-p},\\
&&L_1\cdot(W_p\otimes W_{-p})=(1-p)W_{p+1}\otimes
W_{-p}+(1+p)W_p\otimes W_{1-p}.
\end{eqnarray*}
Denote
\begin{eqnarray*}
&&M_{1}=\max\{\,|p|\,\big|\,a_{1,p}\ne0\},\,\,\ \
M_{2}=\max\{\,|p|\,\big|\,b_{1,p}\ne0\},\\
&&M_{3}=\max\{\,|p|\,\big|\,c_{1,p}\ne0\},\,\,\ \
M_{4}=\max\{\,|p|\,\big|\,d_{1,p}\ne0\}.
\end{eqnarray*}
Using the induction on $\sum_{i=1}^4M_{i}$ in the above identities,
and replacing $D_0$ by $D_0-u_{\rm inn}$, where $u$ is a combination
of some $L_{p}\otimes L_{-p}$, $L_{p}\otimes W_{-p}$, $W_{p}\otimes
L_{-p}$, $W_{p}\otimes W_{-p}$, $L_{1}\otimes c$, $c\otimes L_{1}$,
$W_{p}\otimes c$ and $c\otimes W_{1}$, one can safely suppose
\begin{eqnarray*}
&&a_1=b_1=c_1=d_1=0,\\
&&a_{1,p}=b_{1,p}=c_{1,p}=d_{1,p}=0\ \ \ {\rm if}\ \ p\neq-1,2.
\end{eqnarray*}
Thus the expression of $D_0(L_1)$ can be simplified as (\,recalling
Claim \ref{clai2})
\begin{eqnarray}
D_0(L_1)\!\!\!&\equiv&\!\!\!a_{1,-1}L_{-1}\otimes
L_{2}+a_{1,2}L_2\otimes L_{-1}
+b_{1,-1}L_{-1}\otimes W_{2}+b_{1,2}L_2\otimes W_{-1}\nonumber\\
&&\!\!\!+c_{1,-1}W_{-1}\otimes L_{2}+c_{1,2}W_2\otimes L_{-1}
+d_{1,-1}W_{-1}\otimes W_{2}+d_{1,2}W_2\otimes
W_{-1}.\label{0801261}
\end{eqnarray}
Applying $D_0$ to $[\,L_{-1},L_{1}]=2L_0$, under modulo $\C(c\otimes
c)$, we obtain

\begin{eqnarray*}
&&\mbox{$\sum\limits_{p\in\Z}$}\big((2-p)a_{-1,p-1}+(2+p)a_{-1,p}\big)L_p\otimes
L_{-p}+3a_{1,-1}L_{-1}\otimes L_{1}+3a_{1,2}L_1\otimes
L_{-1}\\
&&+\mbox{$\sum\limits_{p\in\Z}$}\big((2-p)b_{-1,p-1}+(2+p)b_{-1,p}\big)L_p\otimes
W_{-p}+3b_{1,-1}L_{-1}\otimes W_{1}+3b_{1,2}L_1\otimes W_{-1}\nonumber\\
&&+\mbox{$\sum\limits_{p\in\Z}$}\big((2-p)c_{-1,p-1}+(2+p)c_{-1,p}\big)W_p\otimes
L_{-p}+3c_{1,-1}W_{-1}\otimes L_{1}+3c_{1,2}W_1\otimes L_{-1}\\
&&+\mbox{$\sum\limits_{p\in\Z}$}\big((2-p)d_{-1,p-1}+(2+p)d_{-1,p}\big)W_p\otimes
W_{-p}+3d_{1,-1}W_{-1}\otimes W_{1}+3d_{1,2}W_1\otimes
W_{-1}\\
&&+2a_{-1}L_{0}\otimes c+2b_{-1} c\otimes L_{0} +2c_{-1}W_{0}\otimes
c+2d_{-1} c\otimes W_{0}=0.
\end{eqnarray*}
For any $p\in\Z$, comparing the coefficients of $L_p\otimes L_{-p}$,
$L_p\otimes W_{-p}$, $W_p\otimes L_{-p}$, $W_p\otimes W_{-p}$,
$L_{0}\otimes c$, $c\otimes L_{0}$, $W_{0}\otimes c$ and $c\otimes
W_{0}$ respectively in the above equation, on has
\begin{eqnarray*}
&&a_{-1}=b_{-1}=c_{-1}=d_{-1}=0,\\
&&3a_{1,2}+a_{-1,0}+3a_{-1,1}=3b_{1,2}+b_{-1,0}+3b_{-1,1}=0,\\
&&3c_{1,2}+c_{-1,0}+3c_{-1,1}=3d_{1,2}+d_{-1,0}+3d_{-1,1}=0,\\
&&3a_{1,-1}+3a_{-1,-2}+a_{-1,-1}=3b_{1,-1}+3b_{-1,-2}+b_{-1,-1}=0,\\
&&3c_{1,-1}+3c_{-1,-2}+c_{-1,-1}=3d_{1,-1}+3d_{-1,-2}+d_{-1,-1}=0,\\
&&(p-2)a_{-1,p-1}-(p+2)a_{-1,p}=(p-2)b_{-1,p-1}-(p+2)b_{-1,p}=0,\ \
p\neq\pm1,\\
&&(p-2)c_{-1,p-1}-(p+2)c_{-1,p}=(p-2)d_{-1,p-1}-(p+2)d_{-1,p}=0,\ \
p\neq\pm1,
\end{eqnarray*}
which give the following identities:
\begin{eqnarray*}
&&3a_{1,2}+a_{-1,0}+3a_{-1,1}=3b_{1,2}+b_{-1,0}+3b_{-1,1}=0,\\
&&3c_{1,2}+c_{-1,0}+3c_{-1,1}=3d_{1,2}+d_{-1,0}+3d_{-1,1}=0,\\
&&3a_{1,-1}+3a_{-1,-2}-a_{-1,0}=3b_{1,-1}+3b_{-1,-2}-b_{-1,0}=0,\\
&&3c_{1,-1}+3c_{-1,-2}-c_{-1,0}=3d_{1,-1}+3d_{-1,-2}-d_{-1,0}=0,\\
&&a_{-1,p}=b_{-1,p}=c_{-1,p}=d_{-1,p}=0,\ \ \forall\,\,p\in\Z,\,p\neq-2,-1,0,1,\\
&&a_{-1,-1}+a_{-1,0}=b_{-1,-1}+b_{-1,0}=c_{-1,-1}+c_{-1,0}=d_{-1,-1}+d_{-1,0}=0.
\end{eqnarray*}
Thus $D_0(L_{-1})$ and $D_0(L_1)$ can respectively be rewritten as
\begin{eqnarray*}
D_0(L_{-1})\!\!\!&\equiv&\!\!\!a_{-1,-2}L_{-2}\otimes
L_{1}-a_{-1,0}L_{-1}\otimes L_{0}+a_{-1,0}L_0\otimes
L_{-1}+a_{-1,1}L_1\otimes
L_{-2}\\
&&\!\!\!+b_{-1,-2}L_{-2}\otimes W_{1}-b_{-1,0}L_{-1}\otimes W_{0}
+b_{-1,0}L_0\otimes W_{-1}+b_{-1,1}L_1\otimes W_{-2}\\
&&\!\!\!+c_{-1,-2}W_{-2}\otimes L_{1}-c_{-1,0}W_{-1}\otimes L_{0}
+c_{-1,0}W_0\otimes L_{-1}+c_{-1,1}W_1\otimes L_{-2}\\
&&\!\!\!+d_{-1,-2}W_{-2}\otimes W_{1}-d_{-1,0}W_{-1}\otimes W_{0}
+d_{-1,0}W_0\otimes W_{-1}+d_{-1,1}W_1\otimes W_{-2},\\
D_0(L_1)\!\!\!&\equiv&\!\!\!\big(\frac{a_{-1,0}}{3}-a_{-1,-2}\big)L_{-1}\otimes
L_{2}-\big(\frac{a_{-1,0}}{3}+a_{-1,1}\big)L_2\otimes L_{-1}\\
&&\!\!\!+\big(\frac{b_{-1,0}}{3}-b_{-1,-2}\big)L_{-1}\otimes W_{2}
-\big(\frac{b_{-1,0}}{3}+b_{-1,1}\big)L_2\otimes W_{-1}\nonumber\\
&&\!\!\!+\big(\frac{c_{-1,0}}{3}-c_{-1,-2}\big)W_{-1}\otimes
L_{2}-\big(\frac{c_{-1,0}}{3}+c_{-1,1}\big)W_2\otimes L_{-1}\\
&&\!\!\!+\big(\frac{d_{-1,0}}{3}-d_{-1,-2}\big)W_{-1}\otimes
W_{2}-\big(\frac{d_{-1,0}}{3}+d_{-1,1}\big)W_2\otimes W_{-1}.
\end{eqnarray*}
Applying $D_0$ to $[\,L_{2},L_{-1}]=3L_{1}$, under modulo
$\C(c\otimes c)$, we obtain
\begin{eqnarray*}
&&2a_{-1,0}L_2\otimes L_{-1}+3a_{-1,0}L_0\otimes
L_{1}+a_{-1,1}L_3\otimes L_{-2}+4a_{-1,1}L_1\otimes
L_{0}\\
&&+4a_{-1,-2}L_{0}\otimes L_{1}+a_{-1,-2}L_{-2}\otimes
L_{3}-3a_{-1,0}L_{1}\otimes L_{0}-2a_{-1,0}L_{-1}\otimes
L_{2}\\
&&+2b_{-1,0}L_2\otimes W_{-1}+3b_{-1,0}L_0\otimes W_{1}
+b_{-1,1}L_3\otimes W_{-2}+4b_{-1,1}L_1\otimes W_{0}\\
&&+4b_{-1,-2}L_{0}\otimes W_{1}+b_{-1,-2}L_{-2}\otimes W_{3}
-3b_{-1,0}L_{1}\otimes W_{0}-2b_{-1,-1}L_{-1}\otimes W_{2}\\
\end{eqnarray*}
\begin{eqnarray*}
&&+2c_{-1,0}W_2\otimes L_{-1}+3c_{-1,0}W_0\otimes L_{1}
+c_{-1,1}W_3\otimes L_{-2}+4c_{-1,1}W_1\otimes L_{0}\\
&&+4d_{-1,-2}W_{0}\otimes W_{1}+d_{-1,-2}W_{-2}\otimes W_{3}
-3d_{-1,0}W_{1}\otimes W_{0}-2d_{-1,0}W_{-1}\otimes W_{2}\\
&&+4c_{-1,-2}W_{0}\otimes L_{1}+c_{-1,-2}W_{-2}\otimes L_{3}
-3c_{-1,0}W_{1}\otimes L_{0}-2c_{-1,0}W_{-1}\otimes L_{2}\\
&&+2d_{-1,0}W_2\otimes W_{-1}+3d_{-1,0}W_0\otimes W_{1}
+d_{-1,1}W_3\otimes W_{-2}+4d_{-1,1}W_1\otimes W_{0}\\
&&-\big(a_{-1,0}-3a_{-1,-2}\big)L_{-1}\otimes
L_{2}+\big(a_{-1,0}+3a_{-1,1}\big)L_2\otimes
L_{-1}-\big(b_{-1,0}-3b_{-1,-2}\big)L_{-1}\otimes W_{2}
\\
&&+\big(b_{-1,0}+3b_{-1,1}\big)L_2\otimes W_{-1}
-\big(c_{-1,0}-3c_{-1,-2}\big)W_{-1}\otimes
L_{2}+\big(c_{-1,0}+3c_{-1,1}\big)W_2\otimes L_{-1}\\
&&-\big(d_{-1,0}-3d_{-1,-2}\big)W_{-1}\otimes
W_{2}+\big(d_{-1,0}+3d_{-1,1}\big)W_2\otimes
W_{-1}+\mbox{$\sum\limits_{p\in\Z}$}(1+p)a_{2,p}L_{p-1}\otimes
L_{2-p}\\
&&+\mbox{$\sum\limits_{p\in\Z}$}(3-p)a_{2,p}L_p\otimes
L_{1-p}+\mbox{$\sum\limits_{p\in\Z}$}(1+p)b_{2,p}L_{p-1}\otimes
W_{2-p}+\mbox{$\sum\limits_{p\in\Z}$}(3-p)b_{2,p}L_p\otimes
W_{1-p}\\
&&+\mbox{$\sum\limits_{p\in\Z}$}(1+p)c_{2,p}W_{p-1}\otimes
L_{2-p}+\mbox{$\sum\limits_{p\in\Z}$}(3-p)c_{2,p}W_p\otimes
L_{1-p}+\mbox{$\sum\limits_{p\in\Z}$}(1+p)d_{2,p}W_{p-1}\otimes
W_{2-p}\\
&&+\mbox{$\sum\limits_{p\in\Z}$}(3-p)d_{2,p}W_p\otimes
W_{1-p}+3a_2L_{1}\otimes c+3b_2c\otimes L_{1} +3c_2W_{1}\otimes
c+3d_2c\otimes W_{1}=0.
\end{eqnarray*}
For any $p\in\Z$, comparing the coefficients of $L_{1}\otimes c$,
$c\otimes L_{1}$, $W_{1}\otimes c$, $c\otimes W_{1}$, $L_p\otimes
L_{1-p}$, $L_p\otimes W_{1-p}$, $W_p\otimes L_{1-p}$ and $W_p\otimes
W_{1-p}$ respectively in the above equation, we firstly obtain
\begin{eqnarray*}
&&a_2L_{1}\otimes c=b_2c\otimes L_{1}=c_2W_{1}\otimes c=d_2c\otimes W_{1}=0,\\
&&\mbox{$\sum\limits_{p\in\Z}$}(1+p)a_{2,p}L_{p-1}\otimes
L_{2-p}+\mbox{$\sum\limits_{p\in\Z}$}(3-p)a_{2,p}L_p\otimes
L_{1-p}\\
&&+a_{-1,-2}L_{-2}\otimes L_{3}+3(a_{-1,-2}-a_{-1,0})L_{-1}\otimes
L_{2}+(4a_{-1,-2}+3a_{-1,0})L_0\otimes
L_{1}\\
&&+(4a_{-1,1}-3a_{-1,0})L_{1}\otimes
L_{0}+3(a_{-1,0}+a_{-1,1})L_2\otimes L_{-1}+a_{-1,1}L_3\otimes
L_{-2}=0,\\
&&\mbox{$\sum\limits_{p\in\Z}$}(1+p)b_{2,p}L_{p-1}\otimes
W_{2-p}+\mbox{$\sum\limits_{p\in\Z}$}(3-p)b_{2,p}L_p\otimes
W_{1-p}\\
&&+b_{-1,-2}L_{-2}\otimes W_{3}+3(b_{-1,-2}-b_{-1,0})L_{-1}\otimes
W_{2}+(4b_{-1,-2}+3b_{-1,0})L_0\otimes
W_{1}\\
&&+(4b_{-1,1}-3b_{-1,0})L_{1}\otimes
W_{0}+3(b_{-1,0}+b_{-1,1})L_2\otimes W_{-1}+b_{-1,1}L_3\otimes
W_{-2}=0,\\
&&\mbox{$\sum\limits_{p\in\Z}$}(1+p)c_{2,p}W_{p-1}\otimes
L_{2-p}+\mbox{$\sum\limits_{p\in\Z}$}(3-p)c_{2,p}W_p\otimes
L_{1-p}\\
&&+c_{-1,-2}W_{-2}\otimes L_{3}+3(c_{-1,-2}-c_{-1,0})W_{-1}\otimes
L_{2}+(4c_{-1,-2}+3c_{-1,0})W_0\otimes
L_{1}\\
&&+(4c_{-1,1}-3c_{-1,0})W_{1}\otimes
L_{0}+3(c_{-1,0}+c_{-1,1})W_2\otimes L_{-1}+c_{-1,1}W_3\otimes
L_{-2}=0,\\
&&\mbox{$\sum\limits_{p\in\Z}$}(1+p)a_{2,p}W_{p-1}\otimes
W_{2-p}+\mbox{$\sum\limits_{p\in\Z}$}(3-p)a_{2,p}W_p\otimes
W_{1-p}\\
&&+a_{-1,-2}W_{-2}\otimes W_{3}+3(a_{-1,-2}-a_{-1,0})W_{-1}\otimes
W_{2}+(4a_{-1,-2}+3a_{-1,0})W_0\otimes
W_{1}\\
&&+(4a_{-1,1}-3a_{-1,0})W_{1}\otimes
W_{0}+3(a_{-1,0}+a_{-1,1})W_2\otimes W_{-1}+a_{-1,1}W_3\otimes
W_{-2}=0,
\end{eqnarray*}
Then for any $p\in\Z,\,p\neq-2,-1,0,1,2,3$, one has
\begin{eqnarray*}
&&a_2=b_2=c_2=d_2=0,\\
\end{eqnarray*}
\begin{eqnarray*}
&&a_{-1,-2}+5a_{2,-2}=a_{-1,1}+5a_{2,4}
=(p+2)a_{2,p+1}-(p-3)a_{2,p}=0,\\
&&4a_{-1,1}-3a_{-1,0}+3a_{2,2}+2a_{2,1}=3(a_{-1,0}+a_{-1,1})+4a_{2,3}+a_{2,2}=0,\\
&&3a_{-1,-2}-3a_{-1,0}+a_{2,0}+4a_{2,-1}=4a_{-1,-2}+3a_{-1,0}+2a_{2,1}+3a_{2,0}=0,\\
&&b_{-1,-2}+5b_{2,-2}=b_{-1,1}+5b_{2,4}
=(p+2)b_{2,p+1}-(p-3)b_{2,p}=0,\\
&&4b_{-1,1}-3b_{-1,0}+3b_{2,2}+2b_{2,1}=3(b_{-1,0}+b_{-1,1})+4b_{2,3}+b_{2,2}=0,\\
&&3b_{-1,-2}-3b_{-1,0}+b_{2,0}+4b_{2,-1}=4b_{-1,-2}+3b_{-1,0}+2b_{2,1}+3b_{2,0}=0,\\
&&c_{-1,-2}+5c_{2,-2}=c_{-1,1}+5c_{2,4}
=(p+2)c_{2,p+1}-(p-3)c_{2,p}=0,\\
&&4c_{-1,1}-3c_{-1,0}+3c_{2,2}+2c_{2,1}=3(c_{-1,0}+c_{-1,1})+4c_{2,3}+c_{2,2}=0,\\
&&3c_{-1,-2}-3c_{-1,0}+c_{2,0}+4c_{2,-1}=4c_{-1,-2}+3c_{-1,0}+2c_{2,1}+3c_{2,0}=0,\\
&&d_{-1,-2}+5d_{2,-2}=d_{-1,1}+5d_{2,4}
=(p+2)d_{2,p+1}-(p-3)d_{2,p}=0,\\
&&4d_{-1,1}-3d_{-1,0}+3d_{2,2}+2d_{2,1}=3(d_{-1,0}+d_{-1,1})+4d_{2,3}+d_{2,2}=0,\\
&&3d_{-1,-2}-3d_{-1,0}+d_{2,0}+4d_{2,-1}=4d_{-1,-2}+3d_{-1,0}+2d_{2,1}+3d_{2,0}=0,\\
\end{eqnarray*}
which together give the following identities:
\begin{eqnarray*}
&&a_{-1,-2}=-5a_{2,-2}=a_{2,p}=0,\ \
a_{2,-1}=\frac{1}{4}(3a_{-1,0}-a_{2,0}),\ \
a_{2,1}=-\frac{3}{2}(a_{2,0}+a_{-1,0}),\\
&&a_{2,4}=-\frac{1}{5}a_{-1,1},\ \
a_{2,2}=a_{2,0}+2a_{-1,0}-\frac{4}{3}a_{-1,1},\ \
a_{2,3}=-\frac{1}{4}a_{2,0}-\frac{5}{4}a_{-1,0}-\frac{5}{12}a_{-1,1},\\
&&b_{-1,-2}=-5b_{2,-2}=b_{2,p}=0,\ \
b_{2,-1}=\frac{1}{4}(3b_{-1,0}-b_{2,0}),\ \
b_{2,1}=-\frac{3}{2}(b_{2,0}+b_{-1,0}),\\
&&b_{2,4}=-\frac{1}{5}b_{-1,1},\ \
b_{2,2}=b_{2,0}+2b_{-1,0}-\frac{4}{3}b_{-1,1},\ \
b_{2,3}=-\frac{1}{4}b_{2,0}-\frac{5}{4}b_{-1,0}-\frac{5}{12}b_{-1,1},\\
&&c_{-1,-2}=-5c_{2,-2}=c_{2,p}=0,\ \
c_{2,-1}=\frac{1}{4}(3c_{-1,0}-c_{2,0}),\ \
c_{2,1}=-\frac{3}{2}(c_{2,0}+c_{-1,0}),\\
&&c_{2,4}=-\frac{1}{5}c_{-1,1},\ \
c_{2,2}=c_{2,0}+2c_{-1,0}-\frac{4}{3}c_{-1,1},\ \
c_{2,3}=-\frac{1}{4}c_{2,0}-\frac{5}{4}c_{-1,0}-\frac{5}{12}c_{-1,1},\\
&&d_{-1,-2}=-5d_{2,-2}=d_{2,p}=0,\ \
d_{2,-1}=\frac{1}{4}(3d_{-1,0}-d_{2,0}),\ \
d_{2,1}=-\frac{3}{2}(d_{2,0}+d_{-1,0}),\\
&&d_{2,4}=-\frac{1}{5}d_{-1,1},\ \
d_{2,2}=d_{2,0}+2d_{-1,0}-\frac{4}{3}d_{-1,1},\ \
d_{2,3}=-\frac{1}{4}d_{2,0}-\frac{5}{4}d_{-1,0}-\frac{5}{12}d_{-1,1},
\end{eqnarray*}
for any $p\in\Z,\,p\neq-1,0,1,2,3$. Then one can rewrite
$D_0(L_{-1})$, $D_0(L_1)$ and $D_0(L_2)$ respectively as

\begin{eqnarray*}
D_0(L_{-1})\!\!\!&\equiv&\!\!\!a_{-1,0}(L_0\otimes
L_{-1}-L_{-1}\otimes L_{0})+b_{-1,0}(L_0\otimes W_{-1}-L_{-1}\otimes
W_{0})
\\
&&\!\!\!+c_{-1,0}(W_0\otimes L_{-1}-W_{-1}\otimes
L_{0})+d_{-1,0}(W_0\otimes W_{-1}-W_{-1}\otimes W_{0}),\\
\end{eqnarray*}
\begin{eqnarray*}
D_0(L_1)\!\!\!&\equiv&\!\!\!\frac{a_{-1,0}}{3}(L_{-1}\otimes
L_{2}-L_2\otimes L_{-1})+\frac{b_{-1,0}}{3}(L_{-1}\otimes W_{2}
-L_2\otimes W_{-1})\nonumber\\
&&\!\!\!+\frac{c_{-1,0}}{3}(W_{-1}\otimes L_{2}-W_2\otimes
L_{-1})+\frac{d_{-1,0}}{3}(W_{-1}\otimes W_{2}-W_2\otimes W_{-1}),\\
D_0(L_{2})\!\!\!&\equiv&\!\!\!\frac{1}{4}(3a_{-1,0}-a_{2,0})L_{-1}\otimes
L_{3}+a_{2,0}L_0\otimes
L_{2}-\frac{3}{2}(a_{2,0}+a_{-1,0})L_1\otimes
L_{1}\\
&&\!\!\!+(a_{2,0}+2a_{-1,0})L_2\otimes
L_{0}-\frac{1}{4}(a_{2,0}+5a_{-1,0})L_3\otimes
L_{-1}\\
&&\!\!\!+\frac{1}{4}(3b_{-1,0}-b_{2,0})L_{-1}\otimes
W_{3}+b_{2,0}L_0\otimes W_{2}
-\frac{3}{2}(b_{2,0}+b_{-1,0})L_1\otimes W_{1}\\
&&\!\!\!+(b_{2,0}+2b_{-1,0})L_2\otimes W_{0}-\frac{1}{4}(b_{2,0}+5b_{-1,0})L_3\otimes W_{-1}\\
&&\!\!\!+\frac{1}{4}(3c_{-1,0}-c_{2,0})W_{-1}\otimes
L_{3}+c_{2,0}W_0\otimes L_{2}
-\frac{3}{2}(c_{2,0}+c_{-1,0})W_1\otimes L_{1}\\
&&\!\!\!+(c_{2,0}+2c_{-1,0})W_2\otimes L_{0}-\frac{1}{4}(c_{2,0}+5c_{-1,0})W_3\otimes L_{-1}\\
&&\!\!\!+\frac{1}{4}(3d_{-1,0}-d_{2,0})W_{-1}\otimes
W_{3}+d_{2,0}W_0\otimes W_{2}
-\frac{3}{2}(d_{2,0}+d_{-1,0})W_1\otimes W_{1}\\
&&\!\!\!+(d_{2,0}+2d_{-1,0})W_2\otimes
W_{0}-\frac{1}{4}(d_{2,0}+5d_{-1,0})W_3\otimes W_{-1}.
\end{eqnarray*}
Applying $D_0$ to $[\,L_{1},L_{-2}]=3L_{-1}$, under modulo
$\C(c\otimes c)$, we obtain

\begin{eqnarray*}
&&\mbox{$\sum\limits_{p\in\Z}3(1-p)a_{-2,p}$}L_{p+1}\otimes
L_{-2-p}+3\mbox{$\sum\limits_{p\in\Z}(3+p)a_{-2,p}$}L_p\otimes
L_{-1-p}\\
&&+3\mbox{$\sum\limits_{p\in\Z}(1-p)b_{-2,p}$}L_{1+p}\otimes
W_{-2-p}+3\mbox{$\sum\limits_{p\in\Z}(3+p)b_{-2,p}$}L_p\otimes
W_{-1-p}\\
&&+3\mbox{$\sum\limits_{p\in\Z}(1-p)c_{-2,p}$}W_{1+p}\otimes
L_{-2-p}+3\mbox{$\sum\limits_{p\in\Z}(3+p)c_{-2,p}$}W_p\otimes
L_{-1-p}\\
&&+3\mbox{$\sum\limits_{p\in\Z}(1-p)d_{-2,p}$}W_{1+p}\otimes
W_{-2-p}+3\mbox{$\sum\limits_{p\in\Z}(3+p)d_{-2,p}$}W_p\otimes
W_{-1-p}\\
&&+9a_{-2}L_{-1}\otimes c+9b_{-2} c\otimes L_{-1}
+9c_{-2}W_{-1}\otimes c+9d_{-2} c\otimes W_{-1}\\
&&-13a_{-1,0}L_0\otimes L_{-1}-a_{-1,0}L_2\otimes
L_{-3}+a_{-1,0}L_{-3}\otimes L_{2}+13a_{-1,0}L_{-1}\otimes
L_{0}\\
&&-13b_{-1,0}L_0\otimes W_{-1}-b_{-1,0}L_2\otimes W_{-3}
+b_{-1,0}L_{-3}\otimes W_{2}+13b_{-1,0}L_{-1}\otimes W_{0}\nonumber\\
&&-13c_{-1,0}W_0\otimes L_{-1}-c_{-1,0}W_2\otimes L_{-3}
+c_{-1,0}W_{-3}\otimes L_{2}+13c_{-1,0}W_{-1}\otimes L_{0}\\
&&-13d_{-1,0}W_0\otimes W_{-1}-d_{-1,0}W_2\otimes W_{-3}
+d_{-1,0}W_{-3}\otimes W_{2}+13d_{-1,0}W_{-1}\otimes W_{0}\\
&&-9a_{-1,0}L_0\otimes L_{-1}+9a_{-1,0}L_{-1}\otimes
L_{0}-9b_{-1,0}L_0\otimes W_{-1}+9b_{-1,0}L_{-1}\otimes W_{0}
\\
&&-9c_{-1,0}W_0\otimes L_{-1}+9c_{-1,0}W_{-1}\otimes
L_{0}-9d_{-1,0}W_0\otimes W_{-1}+9d_{-1,0}W_{-1}\otimes W_{0}=0.
\end{eqnarray*}
For any $p\in\Z$, comparing the coefficients of $L_p\otimes
L_{-1-p}$, $L_p\otimes W_{-1-p}$, $W_p\otimes L_{-1-p}$, $W_p\otimes
W_{-1-p}$, $L_{-1}\otimes c$, $c\otimes L_{-1}$, $W_{-1}\otimes c$
and $c\otimes W_{-1}$ respectively in the above equation, one has
\begin{eqnarray*}
&&(p-2)a_{-2,p-1}-(p+3)a_{-2,p}=0,\\
&&a_{-2}=b_{-2}=c_{-2}=d_{-2}=9a_{-2,-2}+6a_{-2,-1}+22a_{-1,0}=0,\\
&&6a_{-2,-1}+9a_{-2,0}-22a_{-1,0}=15a_{-2,2}-a_{-1,0}=
15a_{-2,-4}+a_{-1,0}=0,\\
&&(p-2)b_{-2,p-1}-(p+3)b_{-2,p}=9b_{-2,-2}+6b_{-2,-1}+22b_{-1,0}=0,\\
&&6b_{-2,-1}+9b_{-2,0}-22b_{-1,0}=15b_{-2,2}-b_{-1,0}=
15b_{-2,-4}+b_{-1,0}=0,\\
&&(p-2)c_{-2,p-1}-(p+3)c_{-2,p}=9c_{-2,-2}+6c_{-2,-1}+22c_{-1,0}=0,\\
&&6c_{-2,-1}+9c_{-2,0}-22c_{-1,0}=15c_{-2,2}-c_{-1,0}=
15c_{-2,-4}+c_{-1,0}=0,\\
&&(p-2)d_{-2,p-1}-(p+3)d_{-2,p}=9d_{-2,-2}+6d_{-2,-1}+22d_{-1,0}=0,\\
&&6d_{-2,-1}+9d_{-2,0}-22d_{-1,0}=15d_{-2,2}-d_{-1,0}=
15d_{-2,-4}+d_{-1,0}=0,
\end{eqnarray*}
for any $p\in\Z,\,\neq-3,-1,0,2$, which together force
\begin{eqnarray*}
&&a_{-1,0}=a_{-2,p}=b_{-2,p}=c_{-2,p}=d_{-2,p}=0\ \
\forall\,\,p\in\Z,\,p\neq1,
\end{eqnarray*}
Thus $D_0(L_{-1})$, $D_0(L_1)$, $D_0(L_{-1})$ and $D_0(L_1)$ can
respectively be rewritten as
\begin{eqnarray}
D_0(L_{-1})\!\!\!&\equiv&\!\!\!D_0(L_1)\equiv0,\label{lr01}\\
D_0(L_{-2})\!\!\!&\equiv&\!\!\!a_{-2,1}L_1\otimes
L_{-3}+b_{-2,1}L_1\otimes W_{-3}+c_{-2,1}W_1\otimes
L_{-3}+d_{-2,1}W_1\otimes
W_{-3},\nonumber\\
D_0(L_{2})\!\!\!&\equiv&\!\!\!-\frac{1}{4}a_{2,0}L_{-1}\otimes
L_{3}+a_{2,0}L_0\otimes L_{2}-\frac{3}{2}a_{2,0}L_1\otimes
L_{1}+a_{2,0}L_2\otimes L_{0}\nonumber\\
&&\!\!\!-\frac{1}{4}a_{2,0}L_3\otimes
L_{-1}-\frac{1}{4}b_{2,0}L_{-1}\otimes W_{3}+b_{2,0}L_0\otimes W_{2}
-\frac{3}{2}b_{2,0}L_1\otimes W_{1}\nonumber\\
&&\!\!\!+b_{2,0}L_2\otimes W_{0}-\frac{1}{4}b_{2,0}L_3\otimes W_{-1}
-\frac{1}{4}c_{2,0}W_{-1}\otimes L_{3}+c_{2,0}W_0\otimes L_{2}\nonumber\\
&&\!\!\!-\frac{3}{2}c_{2,0}W_1\otimes L_{1}+c_{2,0}W_2\otimes L_{0}
-\frac{1}{4}c_{2,0}W_3\otimes L_{-1}-\frac{1}{4}d_{2,0}W_{-1}\otimes W_{3}\nonumber\\
&&\!\!\!+d_{2,0}W_0\otimes W_{2} -\frac{3}{2}d_{2,0}W_1\otimes
W_{1}+d_{2,0}W_2\otimes W_{0}-\frac{1}{4}d_{2,0}W_3\otimes
W_{-1}.\nonumber
\end{eqnarray}

Applying $D_0$ to $[\,L_{1},L_{-2}]=3L_{-1}$, under modulo
$\C(c\otimes c)$, we obtain
\begin{eqnarray*}
&&a_{2,0}L_{-3}\otimes L_{3}-8a_{2,0}L_{-2}\otimes
L_{2}+39a_{2,0}L_{-1}\otimes
L_{1}-32a_{2,0}L_0\otimes L_0\\
&&+(27a_{2,0}-20a_{-2,1})L_1\otimes L_{-1}-8a_{2,0}L_2\otimes
L_{-2}+(a_{2,0}-4a_{-2,1})L_3\otimes L_{-3}\\
&&b_{2,0}L_{-3}\otimes W_{3}-8b_{2,0}L_{-2}\otimes
W_{2}+39b_{2,0}L_{-1}\otimes
W_{1}-32b_{2,0}L_0\otimes W_0\\
&&+(27b_{2,0}-20b_{-2,1})L_1\otimes W_{-1}-8b_{2,0}L_2\otimes
W_{-2}+(b_{2,0}-4b_{-2,1})L_3\otimes W_{-3}\\
&&c_{2,0}W_{-3}\otimes L_{3}-8c_{2,0}W_{-2}\otimes
L_{2}+39c_{2,0}W_{-1}\otimes
L_{1}-32c_{2,0}W_0\otimes L_0\\
&&+(27c_{2,0}-20c_{-2,1})W_1\otimes L_{-1}-8c_{2,0}W_2\otimes
L_{-2}+(c_{2,0}-4c_{-2,1})W_3\otimes L_{-3}\\
&&d_{2,0}W_{-3}\otimes W_{3}-8d_{2,0}W_{-2}\otimes
W_{2}+39d_{2,0}W_{-1}\otimes
W_{1}-32d_{2,0}W_0\otimes W_0\\
&&+(27d_{2,0}-20d_{-2,1})W_1\otimes W_{-1}-8d_{2,0}W_2\otimes
W_{-2}+(d_{2,0}-4d_{-2,1})W_3\otimes W_{-3}=0,
\end{eqnarray*}
which together force (\,comparing the coefficients of the tensor
products)
\begin{eqnarray*}
a_{2,0}=b_{2,0}=c_{2,0}=d_{2,0}=a_{-2,1}=b_{-2,1}=c_{-2,1}=d_{-2,1}=0.
\end{eqnarray*}
Thus we can deduce
\begin{eqnarray}\label{lr02}
&&D_0(L_{-2})\equiv D_0(L_2)\equiv0.
\end{eqnarray}
Since the Virasoro subalgebra of $\WW$, denoted by $\mathcal
{V}ir:={\rm Span}\{L_n\,|\,n\in\Z\}$ can be generated by the set
$\{L_{-2},\,L_{-1},\,L_1,\,L_2\}$, then by (\ref{lr01}) and
(\ref{lr02}), one has
\begin{eqnarray}\label{lrL}
&&D_0(L_{n})\equiv 0,\ \ \forall\,\,n\in\Z.
\end{eqnarray}

Applying $D_0$ to $[L_{0},[L_{0},W_2]]=2W_2$ and using (\ref{lrL}),
on has
\begin{eqnarray*}
&&L_0\cdot L_{0}\cdot D_0(W_2)=2D_0(W_2).
\end{eqnarray*}
Then using (\ref{0801252}), we obtain

\begin{eqnarray*}
&&\mbox{$\sum\limits_{p\in\Z}2e_{2,p}$}L_p\otimes
L_{2-p}+\mbox{$\sum\limits_{p\in\Z}2f_{2,p}$}L_p\otimes W_{2-p}+2e_2
L_2\otimes c+2f_2 c\otimes L_2\nonumber\\
&&+\mbox{$\sum\limits_{p\in\Z}2g_{2,p}$}W_p\otimes
L_{2-p}+\mbox{$\sum\limits_{p\in\Z}2h_{2,p}$}W_p\otimes W_{2-p}+2g_2
W_2\otimes c+2h_2 c\otimes W_2\\
&&=\mbox{$\sum\limits_{p\in\Z}$}p^2e_{2,p}L_{p}\otimes
L_{2-p}+\mbox{$\sum\limits_{p\in\Z}$}p(2-p)e_{2,p}L_{p}\otimes
L_{2-p}+\mbox{$\sum\limits_{p\in\Z}$}p(2-p)e_{2,p}L_{p}\otimes
L_{2-p}\\
&&+\mbox{$\sum\limits_{p\in\Z}$}(p-2)^2e_{2,p}L_p\otimes
L_{2-p}+\mbox{$\sum\limits_{p\in\Z}$}p^2f_{2,p}L_{p}\otimes
W_{2-p}+\mbox{$\sum\limits_{p\in\Z}$}p(2-p)f_{2,p}L_{p}\otimes
W_{2-p}\\
&&+\mbox{$\sum\limits_{p\in\Z}$}p(2-p)f_{2,p}L_{p}\otimes W_{2-p}
+\mbox{$\sum\limits_{p\in\Z}$}((p-2)^2f_{2,p}L_p\otimes W_{2-p}+4e_n
L_{2}\otimes c+4f_2 c\otimes L_2\nonumber\\
&&+\mbox{$\sum\limits_{p\in\Z}$}p^2g_{2,p}W_{p}\otimes
L_{2-p}+\mbox{$\sum\limits_{p\in\Z}$}p(2-p)g_{2,p}W_{p}\otimes
L_{2-p}+\mbox{$\sum\limits_{p\in\Z}$}p(2-p)g_{2,p}W_{p}\otimes
L_{2-p}\\
&&+\mbox{$\sum\limits_{p\in\Z}$}(p-2)^2g_{2,p}W_p\otimes
L_{2-p}+\mbox{$\sum\limits_{p\in\Z}$}p^2h_{2,p}W_{p}\otimes W_{2-p}
+\mbox{$\sum\limits_{p\in\Z}$}p(2-p)h_{2,p}W_p\otimes W_{2-p}\\
&&+\mbox{$\sum\limits_{p\in\Z}$}p(2-p)h_{2,p}W_{p}\otimes W_{2-p}
+\mbox{$\sum\limits_{p\in\Z}$}(p-2)^2h_{2,p}W_p\otimes W_{2-p}
+4g_2W_2\otimes c+4h_2 c\otimes W_2.
\end{eqnarray*}
Comparing the coefficients of the tensor products, one immediately
has
\begin{eqnarray*}
e_{2,p}=f_{2,p}=g_{2,p}=h_{2,p}=e_2=f_2=g_2=h_2=0,\ \
\forall\,\,p\in\Z,
\end{eqnarray*}
which implies
\begin{eqnarray}\label{lrW2}
&&D_0(W_{2})\equiv 0.
\end{eqnarray}
According to the fact that the algebra $\WW$ is generated by the set
$\{L_{-2},\,L_1,\,L_2,\,W_2\}$, using (\ref{lrL}) and (\ref{lrW2}),
we obtain
\begin{eqnarray}\label{lrWW}
&&D_0(\WW)\equiv 0.
\end{eqnarray}
This proves Claim \ref{clai3}.\QED
\begin{clai}\adddot \label{clai4} \rm
$D_0=0$.
\end{clai}
\vskip4pt
\par By Claims \ref{clai1}, \ref{clai2}\ and\ \ref{clai3}, we have
$D_0(\WW)\subseteq \C(c\otimes c)$. Since $[\WW,\WW]=\WW$, we obtain
$D_0(\WW)\subseteq\WW\cdot(D_0(\WW))=0$. Then Claim \ref{clai4}
follows.\QED

\begin{clai}\adddot \label{clai5} \rm
 For every $D\in{\rm Der}(\WW,\VV)$, (\ref{summable}) is a finite sum. \end{clai}
\vskip4pt
\par By the above claims, we can suppose $D_n=(v_n)_{\rm
inn}$ for some $v_n\in\VV_n$ and $n\in\Z$. If
$\Z'=\{n\in\Z^*\,|\,v_n\ne0\}$ is an infinite set, we obtain
$D(L_0)=\sum_{n\in\Z'\cup\{0\}}L_0\cdot v_n=-\sum_{n\in\Z'}nv_n$ is
an infinite sum, which is not an element in $\VV$, contradicting
with the fact that $D$ is a derivation from $\WW$ to $\VV$. This
proves Claim \ref{clai5} and the proposition.\QED

\begin{lemm}\adddot \label{lemma3ll} Suppose $v\in\VV$ such that
$x\cdot v\in {\rm Im}(1-\tau)$ for all $x\in\WW.$ Then $v\in {\rm
Im}(1-\tau)$.\end{lemm}

\ni{\it Proof.~} First note that $\WW\cdot {\rm Im}(1-\tau)\subset
{\rm Im}(1-\tau).$  We shall prove that after a number of steps in
each of which $v$ is replaced by $v - u$ for some $u\in {\rm
Im}(1-\tau),$ the zero element is obtained and thus proving that
$v\in{\rm Im}(1-\tau)$. Write $v=\sum_{n\in\Z}v_n.$ Obviously,
\begin{eqnarray}\label{eqrx}
v\in {\rm Im}(1-\tau)\ \,\Longleftrightarrow \ \,v_n\in {\rm
Im}(1-\tau) \mbox{ for all } n\in\Z.
\end{eqnarray}
Then $\sum_{n\in\Z}nv_n=L_0\cdot v\in {\rm Im}(1-\tau)$. By
(\ref{eqrx}), $nv_n\in {\rm Im}(1-\tau),$ in particular, $v_{n}\in
{\rm Im}(1-\tau)$ if $n\ne0$. Thus by replacing $v$ by $v-\sum_{n\in
\Z^*}v_n,$ we can suppose $v=v_0\in\VV_0.$ Now we can write
\begin{eqnarray*}
v\!\!\!&=&\!\!\!\mbox{$\sum\limits_{p\in\Z}a_{p}$}L_p\otimes
L_{-p}+\mbox{$\sum\limits_{p\in\Z}b_{p}$}L_p\otimes W_{-p}+a'_0
L_0\otimes c+b'_0 c\otimes L_0\nonumber\\
&&\!\!\!+\mbox{$\sum\limits_{p\in\Z}c_{p}$}W_p\otimes
L_{-p}+\mbox{$\sum\limits_{p\in\Z}d_{p}$}W_p\otimes W_{-p}+c'_0
W_0\otimes c+d'_0 c\otimes W_0,
\end{eqnarray*}
where all the coefficients of the tensor products are complex
numbers and the sums are all finite. Fix the normal total order on
$\Z$ compatible with its additive group structure. Since the
elements of the form $u_{1,p}:=L_p\otimes L_{-p}-L_{-p}\otimes
L_{p}$, $u_{2,p}:=L_p\otimes W_{-p}-W_{-p}\otimes L_{p}$,
$u_{3,p}:=W_p\otimes W_{-p}-W_{-p}\otimes W_{p}$, $u_{1}:=L_0\otimes
c-c\otimes L_{0}$ and $u_{2}:=W_0\otimes c-c\otimes W_{0}$ are all
in $\in {\rm Im}(1-\tau),$ replacing $v$ by $v-u$, where $u$ is a
combination of some $u_{1,p}$, $u_{2,p}$, $u_{3,p}$, $u_{1}$ and
$u_{2}$, we can suppose
\begin{eqnarray}
&&b'_0=c_p=d'_0=0,\ \ \forall\,\,p\in\Z,\label{wpqr1}\\
&&a_{p}\ne 0\ \,\Longrightarrow\ \,p>0\ \mbox{ or }\ p=0,\label{wpqr2}\\
&&d_{p}\ne 0\ \,\Longrightarrow\ \,p>0\ \mbox{ or }\
p=0.\label{wpqr3}
\end{eqnarray}
Then the $v$ can be rewritten as
\begin{eqnarray}\label{sm1}
v=\mbox{$\sum\limits_{p\in\Z_+}a_{p}$}L_p\otimes
L_{-p}+\mbox{$\sum\limits_{p\in\Z}b_{p}$}L_p\otimes
W_{-p}+\mbox{$\sum\limits_{p\in\Z_+}d_{p}$}W_p\otimes W_{-p}+a'_0
L_0\otimes c+c'_0 W_0\otimes c.
\end{eqnarray}
First assume that $a_{p}\ne 0$ for some $p>0$. Choose $q>0$ such
that $q\ne p$. Then we see that the term $L_{p+q}\otimes L_{-p}$
appears in $L_{q}\cdot v,$ but (\ref{wpqr2}) implies that the term
$L_{-p}\otimes L_{p+q}$ does not appear in $L_q\cdot v$, a
contradiction with the fact that $L_q\cdot v\in {\rm Im}(1-\tau)$.
Then one further can suppose $a_p=0,\ \forall\,\,p\in\Z^*$.
Similarly, one also can suppose $d_p=0,\ \forall\,\,p\in\Z^*$.
Therefore, the identity (\ref{sm1}) becomes
\begin{eqnarray}\label{sm2}
v=a_{0}L_0\otimes L_{0}+d_{0}W_0\otimes W_{0}+a'_0 L_0\otimes c+c'_0
W_0\otimes c+\mbox{$\sum\limits_{p\in\Z}$}b_{p}L_p\otimes W_{-p}.
\end{eqnarray}
Finally, we mainly use the fact ${\rm Im}(1-\tau)\subset{\rm
Ker}(1+\tau)$ and our suppose $\WW\cdot v\subset{\rm Im}(1-\tau)$ to
deduce $a_0=d_0=a'_0=c'_0=b_p=0$, $\forall\,\,p\in\Z$. One has the
computation
\begin{eqnarray*}
0\!\!\!&=&\!\!\!(1+\tau)L_1\cdot v\\
\!\!\!&=&\!\!\!2a_{0}(L_1\otimes L_{0}+L_0\otimes
L_{1})+2d_{0}(W_1\otimes W_{0}+W_0\otimes W_{1})
+a'_0(L_1\otimes c+c\otimes L_1)\\
&&\!\!\!+\mbox{$\sum\limits_{p\in\Z}$}\big((2-p)b_{p-1}
+(1+p)b_{p}\big)L_p\otimes
W_{1-p}+\mbox{$\sum\limits_{p\in\Z}$}\big((2-p)b_{p-1}
+(1+p)b_{p}\big)W_{1-p}\otimes L_p\\
&&\!\!\!+c'_0(W_1\otimes c+c\otimes W_1).
\end{eqnarray*}
Then noticing the set $\{p\,|\,b_p\ne0\}$ of finite rank and
comparing the coefficients of the tensor products, one immediately
gets
\begin{eqnarray*}
&&b_0=-2b_{-1}=-2b_1,\\
&&a_0=d_0=a'_0=c'_0=b_p=0,\ \ \forall\,\,p\in\Z,\,p\neq0,\pm1.
\end{eqnarray*}
Then (\ref{sm2}) can be rewritten as
\begin{eqnarray}\label{sm3}
v=b_{1}(L_{-1}\otimes W_{1}-2L_0\otimes W_{0}+L_1\otimes W_{-1}).
\end{eqnarray}
observing the computation
\begin{eqnarray*}
0\!\!\!&=&\!\!\!(1+\tau)L_2\cdot v\\
\!\!\!&=&\!\!\!b_{1}(1+\tau)L_2\cdot (L_{-1}\otimes
W_{1}-2L_0\otimes W_{0}+L_1\otimes W_{-1})\\
\!\!\!&=&\!\!\!b_{1}(1+\tau)(6L_{1}\otimes W_{1}+L_{-1}\otimes
W_{3}-4L_2\otimes W_{0}-4L_0\otimes W_{2}+L_3\otimes W_{-1}),
\end{eqnarray*}
which forces $b_1=0$, then $b_0=-2b_{-1}=-2b_1=0$. Then the lemma
follows.\QED

\vskip10pt \ni{\it Proof of Theorem \ref{main}.} Let $(\WW
,[\cdot,\cdot],\D)$ be a Lie bialgebra structure on $\WW$. By
(\ref{bLie-d}), (\ref{deriv}) and Proposition \ref{proposition},
$\D=\D_r$ is defined by (\ref{D-r}) for some $r\in\WW\otimes\WW$. By
(\ref{cLie-s-s}), ${\rm Im}\,\D\subset{\rm Im}(1-\tau)$. Thus by
Lemma \ref{lemma3ll}, $r\in{\rm Im}(1-\tau)$. Then (\ref{cLie-j-i}),
(\ref{add-c}) and Corollary \ref{coro1} show that $c(r)=0$. Thus
Definition \ref{def2} says that $(\WW ,[\cdot,\cdot],\D)$ is a
triangular coboundary Lie bialgebra.\QED
\end{document}